\documentclass[dvips,aap]{imsart}
\RequirePackage[OT1]{fontenc}
\RequirePackage{amsthm,amsmath,natbib,amssymb,amsfonts,graphicx}
\RequirePackage[colorlinks]{hyperref}
\startlocaldefs
\numberwithin{equation}{section}
\theoremstyle{plain}
\newtheorem{thm}{Theorem}[section]
\newtheorem{lem}[thm]{Lemma}

\newtheorem{prop}[thm]{Proposition}
\theoremstyle{definition}

\theoremstyle{remark}

\newtheorem{rem}[thm]{Remark}

\endlocaldefs

\begin{document}

\begin{frontmatter}
\title{Global regularity and probabilistic schemes for free boundary surfaces of multivariate American derivatives and their Greeks}
\runtitle{Free boundaries and their Greeks }

\begin{aug}
\author{\fnms{J\"org} \snm{Kampen}\thanksref{}\ead[label=e2]{kampen@wias-berlin.de}}
\runauthor{J. Kampen}

\affiliation{Weierstrass Institute}

\address{Weierstrass Institute for Applied Analysis and Stochastics\\
Mohrenstr. 39
\\
10117, Berlin, Germany\\
\phantom{E-mail: kampen@wias-berlin.de\ }\printead*{e2}
}

\end{aug}

\begin{abstract}
In a rather general setting of multivariate diffusion market models we derive global iterative probabilistic schemes for computing the free boundary and its Greeks for a generic class of American derivative models using front-fixing methods. The convergence of the scheme is closely linked to a proof of global regularity of the free boundary surface. 
\end{abstract}



\end{frontmatter}

\section{Introduction}
Probabilistic representations of solutions are popular in finance because of the curse of dimensionality  which requires the application of Monte-Carlo methods. Especially, this is true for more intricate problems such as the free boundary problems related to pricing, hedging, and the optimal exercise strategies of American options. The present work contributes to these problems as it establishes a stable iterative functional solution scheme in a general framework of diffusion type market models. However, it is also a contribution to the theoretical question of regularity of the free boundary itself where we establish results which are neither covered by the general results of \cite{Ca, CPS} (cf. discusson below) nor by more special recent regularity results on American derivatives in \cite{PS, Sh}.
There are several methods in order to investigate regularity used in the literature. One is to differentiate equations in order to get equations for derivatives. Another is to find estimates for tangent paraboloids of the solutions (cf. \cite{WaI}). A second method is to find the solution in explicit form which is possible in a very limited class of problems only. In \cite{CPS} the multivariate case of the parabolic potential problem is investigated, i.e. the case of constant coefficients and with no assumption on the sign of the solution. The points of the free boundary are classified into regular and singular points by an energy and a density criterion and $C^{\infty}$ regularity is proved around the regular points. Singular points are not considered. Hence, the regularity results of the present work which establishes overall regularity especially for the American basket Put option are not covered even in the case of the multivariate Black-Scholes model. However our results extend to a considerable class of models with variable coefficients. The only condition (beside regularity conditions on the coefficients) is a global graph condition which allows for a global transformation of the free boundary (also known as front-fixing).   
The method of front fixing itself is discussed in several papers on univariate American options (cf. \cite{WK, NST}) where it is less efficient than a moving boundary approach (cf. \cite{Mu}).
We shall investigate the American basket put option since it is the most prominent and important example. Regularity results for the American put option in the univariate case go back to \cite{McK} and \cite{vM} in the case of constant coefficients. In \cite{Fr2} the univariate case is studied in the case of the related Stefan problem, and where coefficients depend on time and space. A more recent contribution for American Options in the univariate case with variable coefficients is given in \cite{BDM}. In the present work we obtain global regularity results for the American basket put.
However, our method can be applied to the other popular options (we shall provide a list below). We transform the free boundary problem to a nonlinear problem on a domain homeomorphic to a cube. Then the solution of the free boundary surface is presented in terms of a nonlinear integral equation involving convolutions with transition densities of linear parabolic equations and linear Volterra equations. Analysis of these equations leads to regularity results which are simple (use of Banach fixed point theorem) and charming. Our main result is the global regularity of the free boundary surface except for the final time $T$ if a global graph condition is satisfied (which can be justified for a considerable class of models). The nonlinear integral equation is also the basis for our probabilistic scheme for computing the free boundary function, its time derivative and its spatial derivatives up to second order (Greeks). The only other Monte-Carlo method for computing the free boundary (not Greeks) of a multivariate American Put option known to the author is described in \cite{MRS,BM}. We reduce the solution of this integral equation to an iteration of the solution of linear equations (essentially two parabolic equations and a linear Volterra integral equation). The corresponding convolutions with the transition density can be computed by WKB approximations (cf. \cite{Ka1,Ka2,Ka3}). Recent developments in the computations of Greeks with Monte-Carlo Methods can be found in  \cite{EFT, FrKa, KKS}. In \cite{KKS} it is shown that WKB approximations can be a very efficient tool for the computation of Greeks for high-dimensional models. The convergence of the resulting iteration method is based on the regularity results for the free boundary.

The outline of the present article is as follows. In Section 2 we describe the general frame work of market models considered. We restrict ourselves to diffusion models (a class of stochastic volatility models can be included). We formulate the free boundary problem to be solved. In principle the analysis could be extended to some class of models based on Levy processes. However the treatment of global operators in the context of free boundaries leads to additional difficulties which will be considered elsewhere. In Section 3 we apply the front fixing method to the market models set up in the previous Section. The free boundary function appears in the coefficients of a nonlinear parabolic operator and is coupled to the boundary conditions. The smooth fit condition of the original free boundary transforms into a mixed boundary condition on a hyperplane. In Section 4 we derive a nonlinear integral equation which characterizes the free boundary and is based on the nonlinear initial boundary value problem characterizing the value function and established in the previous Section. In Section 5 a proof of the global existence and the regularity of the free boundary function is established. This is done via the nonlinear integral equation established in the previous Section. The integral equation is shown to determine a map in the Banach space of H\"older continuous functions with H\"older continuous derivatives up to the second order. Note that interior norms are used because the problem is not differentiable at maturity time and only $C^1$ on the boundary hyperplane where the mixed boundary condition is related to the smooth fit condition of the original problem. The solution of the free boundary function is the fixed point of a nonlinear integral equation in a Banach space. In Section 6 an iteration scheme for the computation of the free boundary function and its Greeks is given. The description is distributed over four subsections. In the first subsection subproblems related to the nonlinear integral equation characterizing the free boundary function are identified. The second subsection recalls some recent results of the WKB-expansion of the transition density (fundamental solution of linear parabolic equations). The third subsection describes the algorithm and the fourth subsection provides a convergence and qualitative error estimate analysis. A conclusion is given in Section 7 where we also indicate future research.

\section{General framework}
Let $O$ be an open domain in or equal to ${\mathbb R}^n$ and let $T\in (0,\infty)$ be the time horizon. We write
$D=(0,T)\times O$ for the domain where $O$ has to be specified case by case. For each starting point $(t,x)\in [0,T]\times O$, consider the following stochastic differential equation 
\begin{equation}\label{stocheq}
\begin{array}{ll}
X_t^{t,x}=x\in O,\\
\\
dX_s^ {t,x,k}=\mu_k(s, X_s^ {t,x})ds+\sum_{j=1}^m\sigma_{kj}(s,X_s^ {t,x})dW^j_s,
\end{array}
\end{equation}
for continuous  drift functions $\mu_k$ and local volatility functions $\sigma_{ij}$ 
\begin{equation}
\mu_k, \sigma_{ij} :[0,T]\times {\mathbb R}^n\rightarrow {\mathbb R},~~i,j,k=1,\cdots ,n,
\end{equation}
with an ${\mathbb R}^m$-valued Brownian motion $W=(W^ j)_{j=1,\cdots,m}$ and where 
$X_s^ {t,x}=(X_s^ {t,x,k})_{k=1,\cdots ,n}$.
We assume  that the functions $\mu_k$ and $\sigma_{jk}$ are locally Lipschitz-continuous in $x$, uniformly in $t$, i.e. for each compact subset $K\subset {\mathbb R}^n$ there is a constant $c$ (dependent on $K$) such that 
\begin{equation}\label{lip}
|\mu_k(t,x)-\mu_k(t,y)|\leq c|x-y|, \mbox{and}~~|\sigma_{ij}(t,x)-\sigma_{ij}(t,y)|\leq c|x-y|
\end{equation}
for all $t\in[0,T]$ $x,y\in \Omega$. The latter assumption implies that \eqref{stocheq} has a unique strong solution for any given filtered probability space $\left( \Omega, {\cal F},({\cal F}_t)_{t\in [0,T]}, P\right) $ and Brownian motion $W$ up to a possibly finite random explosion time. Therefore we add the assumption that $P\left(\sup_{s\in[t,T]}|X^{t,x}_s|<\infty\right)=1.$
Then theorem V. 38 of \cite{Pr} implies that \eqref{stocheq} has a strong solution. In this context we consider a stochastic volatility market model system with a process $S_t$ modeling the price of $n$ assets and an $d$-dimensional background process $Y$ driving the local volatility. Let ${\mathbb R}^m_+$ denote the set of $m$-tuples of strictly positive real numbers. Then for each starting point $(t,x)\in [0,T]\times \tilde{O}\subseteq {\mathbb R}^m_+\times {\mathbb R}^d$ the stochastic differential equation 
\begin{equation}\label{marketincomp}
\begin{array}{ll}
(S_t^{t,x,y},Y_t^ {t,y})=(x,y)\in \tilde{O},\\
\\
\frac{dS_s^ {t,x,k}}{S_s^ {t,x,k}}=r(s, X_s^ {t,x})ds+\sum_{j=1}^n\sigma_{kj}(s,X_s^ {t,x},Y_s^{t,x})dW^j_s,\\
\\
dY_s^{t,y,l}=\nu_l(s, Y_s^ {t,x})ds+\sum_{j=1}^p\beta_{lj}(s,Y_s^{t,x})dW^j_s
\end{array}
\end{equation} 
has a unique strong solution, if the drift functions
\begin{equation}
\begin{array}{ll}
r :[0,T]\times {\mathbb R}^m\rightarrow {\mathbb R},\\
\\
\nu_l :[0,T]\times {\mathbb R}^d\rightarrow {\mathbb R},~~l=1,\cdots ,d,
\end{array}
\end{equation}
and the (volatility of) volatility functions  
\begin{equation}
\beta_{lj} :[0,T]\times {\mathbb R}^d\rightarrow {\mathbb R},~~l\in\{1,\cdots ,p\},~~j\in\{1,\cdots ,d\}
\end{equation}
satisfy conditions of form \eqref{lip}. In this context the price $V(t,x,y)$ of an American option with payoff $\phi$ is given by
\begin{equation}\label{price}
V(t,x,y):=\sup_{\tau \in \mbox{Stop}_{[t,T]}} E_Q\left[\exp\left( -\int_t^ {\tau}r(s, S_t^{t,x,y})\right)\phi\left( \tau , S_{\tau}^{t,x,y} \right)  \right],
\end{equation}
where $\mbox{Stop}_{[t,T]}$ is the set of all ${\cal F}_t$ stopping times with value in $[t,T]$. If the diffusion matrix of the process $Z^{t,x,y}:=(X_t^{t,x},Y_t^ {t,y})$ satisfies a strict
uniform ellipticity condition with some constants $0<\lambda <\Lambda <\infty$, then the analysis below applies in this general framework.
\begin{rem}
Some typical stochastic volatility models such as multivariate extensions of Heston's model (or Wishart-type models) do not satisfy a strict ellipticity condition. There are also semi-elliptic diffusion models where the number of Brownian motions is less than the dimensionality of the problem. Some (but not all) models of this type satisfy the H\"{o}rmander condition. Note that for some class of semi-elliptic equations recent techniques in (\cite{FK}) may be combined with techniques considered here.
\end{rem}
  
The measure $Q$ in (\ref{price}) is some equivalent measure to be chosen if the market is incomplete. It is uniquely determined if the market is complete. If the asset process $S$ does not depend on a background process $Y$, then the market process is complete, especially $S_{s}^{t,x,y}\equiv S_{s}^{t,x}$. For simplicity of notation we deal with this situation in the following, where we note that everything can be generalized to the more general incomplete market model \eqref{marketincomp} without problems if a uniform strict ellipticity condition is satisfied (the choice of the equivalent measure is a problem of its own).
In case some ellipticity condition holds for $\sigma\sigma^ T$ standard stochastic control theory can now be used to show that the value function $(t,x)\rightarrow V(t,x)$ satisfies the nonlinear Cauchy problem
\begin{equation}\label{nlincauchy1}
\left\lbrace\begin{array}{ll}
\max\left\lbrace \frac{\partial u}{\partial t}+Lu, \phi -u \right\rbrace=0, \mbox{ in } [0,T)\times {\mathbb R}^n_+\\
\\
u(T,x)=\phi(T,x), \mbox{ in } \left\lbrace T\right\rbrace \times {\mathbb R}^n_+
\end{array}\right.
\end{equation}
in the viscosity sense, where $v_{ij}=( \sigma\sigma^T)_{ij}$ and where
\begin{equation}
Lu\equiv \frac{1}{2}\sum_{ij}v_{ij}S_iS_j\frac{\partial^ 2 u}{\partial S_i\partial S_j}+r\left(\sum_i S_i\frac{\partial u}{\partial S_i}-u \right)
\end{equation}
 Global solutions for the latter type of equations in the viscosity sense were studied in \cite{Am} even in a far more general context.
It is clear that \eqref{nlincauchy1} formally implies that
\begin{equation}\label{nlincauchy2}
\left\lbrace\begin{array}{ll}
 \frac{\partial u}{\partial t}+Lu \leq 0, \\
\\
\left( \frac{\partial u}{\partial t}+Lu\right) \left(\phi-u \right)=0,\\
\\
 \phi -u \geq 0, \mbox{ in } [0,T)\times {\mathbb R}^n_+ \mbox{ a.e.}\\
\\
u(T,S)=\phi(T,S) \mbox{ on } \left\lbrace T\right\rbrace \times {\mathbb R}^n_+,
\end{array}\right.
\end{equation} 
and where global solutions can be obtained under rather mild conditions by variational methods (cf. \cite{Fried} ). The equivalence of this system of variational inequalities with optimal stopping problems under relatively mild conditions have also be studied extensively (cf. \cite{FlSo} ). The set of points $(t,x)$ where the value function $u$ equals the obstacle $\phi$ is called the exercise region. The complement of the exercise region in the domain $D$ is called the continuation region. In \cite{EK} it is proved in a rather general context that the optimal stopping time is the first time where the value process hits the boundary of the exercise region. Therefore the free boundary surface (or the boundary of the exercise region) is crucial information in dealing with American Option. Especially in the context of hedging information concerning the Greeks near the free boundary is crucial. In the present article we set up a method to obtain the Greeks in a general class of multivariate models. Information of the value function in the continuation region and its Greeks are then a by-product of our method. 

Strict ellipticity of the diffusion matrix function $(v_{ij})$ and H\"older continuity of the functions $v_{ij}$ and $r$ together with bounded Lipschitz continuous data $\phi$ is enough to ensure existence of global solutions in viscosity sense for \eqref{nlincauchy1} (cf. \cite{Am} for global existence  for more general conditions and classes of problems), or, with a bit more regularity assumptions (derivatives in $L^2$), global existence for \eqref{nlincauchy2} can be ensured by classical variational methods in Sobolev spaces. In the following we shall assume that at least the former conditions hold. However in the design of the scheme we shall essentially assume that coefficient functions are $C^{\infty}$ with exponential bounds of the derivatives. More precisely, we shall assume that for two positive constants $0<\lambda<\Lambda<\infty$, and all $(t,x)$ and all $n$-dimensional vectors $\xi\neq 0$ with real entries $\xi_i$ we have $0<\lambda |\xi |^2 \leq \sum_{i,j=1}^nv_{ij}(t,x) \xi_i\xi_j\leq \Lambda |\xi|^2$. Furthermore, we assume that the functions $v_{ij},~1\leq i,j\leq n,~r$ are H\"older continuous (with exponent $\alpha \in (0,1)$) and bounded. This ensures also the global existence of the fundamental solution (transition density).
  
\section{Frontfixing with special focus on American Basket Put options}

The front-fixing method presented here works in principle for a large class of options, e.g.
\begin{itemize}

\item the minimum put with payoff $\left(K-\min\left\lbrace S_1,S_2 \right\rbrace \right)^+$,
\item the spread with payoff $\left(K- (S_1-S_2)  \right)^+$,
\item the put on an index or index spreads $\left(K-\sum_{i=1 }^ n\alpha_i S_i\right)^+$, and
\item similar call options in standard models and other standard options.

\end{itemize}
Since suitable front-fixing depends on the pay-off we shall restrict the analysis to the most popular example which is the American Put on an index. We start with the operator
\begin{equation}\label{orig}
\frac{\partial u}{\partial t}+\frac{1}{2}\sum_{ij}v_{ij}S_iS_j\frac{\partial^ 2 u}{\partial S_i\partial S_j}+r\left(\sum_i S_i\frac{\partial u}{\partial S_i}-u \right), 
\end{equation}
where $v_{ij}=(\sigma\sigma^T)_{ij}$ and $r$ may depend on time $t$ and spatial variables $S$.
Let ${\cal E}\subset [0,T]\times {\mathbb R}^n_+$ denote the exercise region and for each $t\in [0,T]$ let ${\cal E}_t$ denote the $t$-section of the exercise region, i.e. ${\cal E}_t:=\left\lbrace x|(t,S)\in {\cal E}\right\rbrace$. 
In general beside basic standard assumptions on diffusion market models introduced in the previous Section we shall assume that 

\begin{itemize}
\item[(GG)]
for each $t\in [0,T]$ we assume that
$0\in {\cal E}_t$ and that ${\cal E}_t$ is star-shaped with respect to $0$, i.e. for all $S\in {\cal E}_t$ and all $\lambda\in [0,1]$ we assume that $\lambda S\in {\cal E}_t$. 
\begin{rem}
Note that this means that for a fixed "angle" at $S=(S_1,\cdots ,S_n)$, i.e. at
$$
\phi_S:=\left( \frac{S_2}{\sum_{i=1}^n S_i},\cdots ,\frac{S_n}{\sum_{i=1}^n S_i}\right) 
$$
we have one intersection point of the free boundary of the section ${\cal E}_t$ and  the ray through $0$ which is determined by the angle $\phi_S$.
\end{rem}
\end{itemize}
\begin{rem}
We call (GG) the global graph condition. The condition (GG) holds if $x\rightarrow u(t,x)$ is convex, where $(t,x)\rightarrow u(t,x)$ denotes the value function of an American Put. This is a sufficient (not necessary) condition for (GG) to hold. Especially, this condition is satisfied for the multivariate Black-Scholes model (Consider the Snell envelope definition in order to verify convexity).
\end{rem}
Hence, the free boundary can be written in terms of the angles in form
\begin{equation}
(t,\phi_S)\rightarrow F(t,\phi_S). 
\end{equation}
We consider the transformation
\begin{equation}
\begin{array}{ll}
\psi: (0,T)\times {\mathbb R}^ n_+\rightarrow (0,T)\times[1,\infty)\times \left(0,1 \right)^ {n-1},\\
\\
\psi(t,S_1,\cdots ,S_n)= \left( t ,\frac{\sum_{i=1}^n S_i}{F},\frac{S_2}{\sum_{i=1}^n S_i},\cdots ,\frac{S_n}{\sum_{i=1}^n S_i}\right).
\end{array}
\end{equation}
Note that the spatial part of $\psi\left( (0,T)\times {\mathbb R}^ n_+\right) $ is homeomorph to the half space $H_{\geq 1}=\left\lbrace x\in {\mathbb R}^ n|x_1\geq 1 \right\rbrace $. In the following the domain $D$ is the interior of the image of $\psi$, i.e. $D:=(0,T)\times (1,\infty)\times \left(0,1 \right)^ {n-1}$. We have
\begin{equation}
S_1=x_1F\left(1-\sum_{j\geq 2}x_j\right),~S_j=x_jx_1F.
\end{equation}
We get
\begin{equation}\label{frontbd}
\left\lbrace \begin{array}{ll}
u_{t}=\frac{F_t}{F}x_1\frac{\partial u}{\partial x_1}+\frac{1}{2}\sum_{ij}a^F_{ij}\frac{\partial^ 2 u}{\partial x_i\partial x_j}+\sum_j b^F_j\frac{\partial u}{\partial x_j}+r\left( x_1\frac{\partial u}{\partial x_1}-u\right) ,\\
\\
(BC1)~~u(0,\infty,x_2,\cdots ,x_n)=0 \mbox{ on $x_1=\infty$}\\
\\
(BC2)~~u_{x_1}(t,1,x_2,\cdots ,x_n)-u(t,1,x_2,\cdots ,x_n)=-K \\
\hspace{6cm}\mbox{ on $x_1=1$}\\
\\
(BC3)~~F(t,x_2,\cdots ,x_n)=K-u(t,1,x_2,\cdots,x_n)\\
\\
(IC)~~u(0,x)=\max \{K-x_1, 0\}
\end{array}\right.
\end{equation}
\begin{rem}
We include (BC1) as an implicit boundary condition in order to indicate that \eqref{frontbd} is equivalent to an initial-boundary value problem of the second type on a finite domain (just by suitable additional transformation with respect to the variable $x_1$).
\end{rem}
The mixed condition $(BC2)$ follows from the smooth fit condition together with $(BC3)$ .
\begin{rem}
Note that in the context of market models based on Levy processes or, more generally, Feller processes the smooth fit condition does not hold in general and one has to be careful concerning generalization at this point. 
\end{rem}
In order to determine the coefficients $a^F_{ij}$ and $b^F_i$ we compute first
\begin{equation}
\begin{array}{ll}
F\frac{\partial x_j}{\partial S_i}=\frac{\delta_{ij}-x_j}{x_1},~
F\frac{\partial x_j}{\partial S_1}=1-\sum_{j\geq 2}(\delta_{ij}-x_j)\frac{F_j}{F},~
\frac{\partial}{\partial S_i}=\sum_j \frac{\partial x_j}{\partial S_i}\frac{\partial}{\partial x_j}.
\end{array}
\end{equation}
We observe that
\begin{equation}
\sum_i S_i\frac{\partial x_j}{\partial S_i}=\sum_i S_i\frac{\delta_{ij}-x_j}{x_1F}=0.
\end{equation}
It follows that
\begin{equation}
\begin{array}{ll}
\sum_i S_i\frac{\partial }{\partial S_i}&=\sum_{ij}\frac{\partial x_j}{\partial S_i}\frac{\partial }{\partial S_i}¬=\sum_{i}S_i\frac{\partial x_1}{\partial S_i}\frac{\partial }{\partial x_1}    +\sum_{j\geq 2}\left(\sum_i S_i\frac{\partial x_1}{\partial S_i} \right) \frac{\partial }{\partial x_j}     \\
\\
&=\sum_i S_i\frac{\partial x_1}{\partial S_i}\frac{\partial}{\partial x_1}=\sum_i S_i\left(\frac{1}{F}-\sum_{j\geq 2}\left(\delta_{ij}-x_j\right)\frac{F_j}{F^ 2} \right) \frac{\partial}{\partial x_1}\\
\\
&=x_1\frac{\partial }{\partial x_1}-\left(\sum_{i}S_i\right) \left(\sum_{j\geq 2}\left(\frac{\partial x_j}{\partial S_i}\right)\frac{F_j}{F}x_1 \right)\frac{\partial }{\partial x_1}=x_1\frac{\partial }{\partial x_1}.
\end{array}
\end{equation}
Hence, we have
\begin{equation}
r\left(\sum_i S_i \frac{\partial}{\partial S_i}\right)=rx_1\frac{\partial}{\partial x_1}.
\end{equation}
It is clear that
\begin{equation}
a^F_{ij}=\sum_{kl}v_{kl}S_kS_l\frac{\partial x_i}{\partial S_l}\frac{\partial x_l}{\partial S_k},~~
b^F_{j}=\sum_{kl}v_{kl}S_kS_l\frac{\partial^2 x_j}{\partial S_k\partial S_l}.
\end{equation}
In order to determine the latter coefficient functions we compute
\begin{equation}
\frac{\partial x_j}{\partial S_i}=\frac{1}{F}\frac{\delta_{ij}-x_j}{x_1}, j\geq 2,~~
\frac{\partial x_1}{\partial S_i}=\frac{1}{F}\left( 1-\sum_{j\geq 2}(\delta_{ij}-x_j)\frac{F_j}{F}\right).
\end{equation}
Next, for $j\geq 2$ we have
\begin{equation}
\frac{\partial^ 2 x_j}{\partial S_i\partial S_k}=\sum_l\frac{\partial \left(\frac{\delta_{ij}-x_j}{Fx_1} \right) }{\partial x_l}, \mbox{ and }
\end{equation}
\begin{equation}
\frac{\partial \left(\frac{\delta_{ij}-x_j}{Fx_1} \right)}{\partial x_l}
=-\frac{\delta_{jl}}{x_1F}+\frac{(x_j-\delta_{ij})(\delta_{1l}F+x_1F_l(1-\delta_{1l}))}{(x_1F)^ 2}.
\end{equation}
Finally,
\begin{equation}
\frac{\partial^ 2 x_1}{\partial S_i\partial S_k}=\sum_l \frac{\partial}{\partial x_l}\left(\frac{1}{F}-\sum_{j\geq 2}(\delta_{ij}-x_j)\frac{F_j}{F^ 2}\right) \frac{\partial x_l}{\partial S_k}, \mbox{ where}
\end{equation}
\begin{equation}
\begin{array}{ll}
\frac{\partial}{\partial x_l}\left(\frac{1}{F}-\sum_{j\geq 2}(\delta_{ij}-x_j)\frac{F_j}{F^ 2}\right)\\
\\
=-\frac{F_l}{F^2}(1-\delta_{1l})-\sum_{j\geq 2}\frac{-\delta_{jl}F_j+(\delta_{ij}-x_j)F_{jl}(1-\delta_{1l})-2F_jF_l(1-\delta_{1l})(\delta_{ij}-x_j)}{F^ 3}.
\end{array}
\end{equation}
Here $\delta_{ij}$ is always the Kronecker Delta, $F_j$ is short for $\frac{\partial F}{\partial x_j}$, $F_{jl}$ is short for $\frac{\partial^ 2 F}{\partial x_j\partial x_l}$. Now we have determined the explicit form of \eqref{frontbd}. The next step is to construct a representation of the solution of \eqref{frontbd} in terms of convolutions with the transition density, i.e the fundamental solution related to \eqref{frontbd}.

\section{Derivation of a system of integral equations characterizing the free boundary}

For a fixed free boundary function $F$ we may consider \eqref{frontbd} as a standard linear initial value boundary problem of the second type and of the form 
\begin{align}\label{init2a}
\frac{\partial u}{\partial t}+\frac{1}{2}\sum_{i,j=1}^{n}a_{ij}\frac
{\partial^{2}u}{\partial x^{i}\partial x^{j}}+\sum_{i=1}^{n}%
b_i\frac{\partial u}{\partial x^{i}}  &  =f, \mbox{ in $D_0$},\\
u(0,x)  &  =g(x) \mbox{ on $\overline{O}$}\\
\frac{\partial}{\partial \nu}u(t,x)+\alpha u(t,x) & =h \mbox{ on $H$},\label{BKa} 
\end{align}
where $H:=\left\lbrace (t,x)\in \overline{D_0}|x_1= 1 \right\rbrace$ is part of the boundary of $D$ ($H$ is for 'hyperplane'), and $D_0=O\times (t_0,T)$ for some $t_0\geq 0$. 
This is an initial value boundary problem, where $\nu$ denotes the inward normal, and $a_{ij}, b_i,f,\alpha,h$ are functions which may depend on time $t$ and the spatial variables $x$. In our case the derivative with respect to the inward normal
reduces to the partial derivative $\frac{\partial}{\partial x_1}$, the function $\alpha$ is the constant function $\alpha\equiv -1$, $h$ is the constant function $h\equiv -K$, and $f\equiv 0$. Furthermore, in case $t_0=0$ the initial condition for the basket Put is just $g(x)=\psi_0(x):=\max\left\lbrace K-x_1,0\right\rbrace $ in the transformed coordinates. In general for $t_0>0$ the initial condition is $\psi_{t_0}(x)=u(t_0,x)$, where $(t,x)\rightarrow u(t,x)$ is the solution of (\ref{frontbd}) on the time interval $\left[t_0,T\right] $. Hence, the equations of \eqref{init2a} simplify to

\begin{align}\label{init2b}
\frac{\partial u}{\partial t}+\frac{1}{2}\sum_{i,j=1}^{n}a_{ij}\frac
{\partial^{2}u}{\partial x^{i}\partial x^{j}}+\sum_{i=1}^{n}%
b_i\frac{\partial u}{\partial x^{i}}   & =0, \mbox{ in $D_0$},\\
u(0,x)    = \max\left\lbrace K-x_1,0\right\rbrace &\mbox{ on $\overline{O}$},\\
\frac{\partial}{\partial x_1}u(t,x)- u(t,x) & =-K \mbox{ on $H$}.\label{BKb} 
\end{align}

We represent the solution in terms of convolutions of the fundamental solution and an integral equation of Volterra type related to the mixed type boundary condition.
First (given a fixed $F$) we define the fundmental solution $p_F\equiv p_F(t,x;\tau,y)$ for $t>\tau$ to solve for each $(\tau,y)$
\begin{align}
\frac{\partial u}{\partial t}+\frac{1}{2}\sum_{i,j=1}^{n}a^F_{ij}(t,x)\frac
{\partial^{2}u}{\partial x^{i}\partial x^{j}}+\sum_{i=1}^{n}%
b_i(x)\frac{\partial u}{\partial x^{i}}  &  =0, \mbox{ in $O\times (0,T]$}\\
u(0,x)  &  =\delta_y(x) \mbox{ at $t=t_0$}.
\end{align}
(where $b_i(x):=b_i^ F+\frac{F_t}{F}x_1+rx_1$). 
The general ansatz for \eqref{init2a}-\eqref{BKa} then is (recall that $\hat{x}_1=(x_2,\cdots,x_n), \hat{y}_1=(y_2,\cdots,y_n)$, and
$dH_y=dy_2dy_3\cdots dy_n$), (writing $H=H_0\times (t_0,T)$) is
 \begin{equation}\label{sol}
\begin{array}{ll}
 u(t,x)=&\int_{t_0}^t \int_{H_0} p_F(t,1,\hat{x}_1;\tau,1,\hat{y}_1)\phi (\tau,1,\hat{y}_1)dH_yd\tau\\
\\
&+\int_O p_F(t,x;0,y)\psi_0(y) dy, 
\end{array}
 \end{equation}
where $\hat{x}_1=(x_2,\cdots,x_n), \hat{y}_1=(y_2,\cdots,y_n)$,
$dH_y=dy_2dy_3\cdots dy_n$). Moreover, we write $H=H_0\times (t_0,T)$.
The boundary condition $\frac{\partial}{\partial x_1}u(t,x)-u(t,x)=-K$ on $H$ ( equal to $\frac{\partial}{\partial x_1}u(t,1,\hat{x}_1)-u(t,1,\hat{x}_1)=-K$) then reduces to the Volterra type equation 
\begin{equation}\label{vol}
\begin{array}{ll}
\frac{1}{2}\phi(t,1,\hat{x}_1)=\Gamma(t,1,\hat{x}_1)+\\
\\ \int_{t_0}^ t \int_{H_0} (\frac{\partial}{\partial x_1}p_F(t,1,\hat{x}_1,\tau,\hat{y}_1)-p_F(t,1,\hat{x}_1;\tau,1,\hat{y}_1))\phi (\tau,1,\hat{y}_1)dH_yd\tau,
\end{array}
\end{equation}
(defining $\phi$) with 
\begin{equation}
\begin{array}{rr}
\Gamma (t,x)=&\int_{O} \left( \frac{\partial}{\partial x_1}p_F(t,x;t_0,y)- p_F(t,x;t_0,y)\right) \psi_0(y)dy+K.
\end{array}
\end{equation}
We apply this to the initial value boundary problem \eqref{frontbd} in order to obtain a system of integral equations for the free boundary. First, from the boundary condition (BC3) in \eqref{frontbd} we get the integral equation
\begin{equation}
\begin{array}{lll}
F(t,x_2,\cdots ,x_n)=&K-u(t,1,\hat{x}_1)=\\
\\
\hspace{2cm}&K-\int_{t_0}^ t \int_{H_0} p(t,1,\hat{x}_1;\tau,1,\hat{y}_1)\phi (\tau,1,\hat{y}_1)dH_yd\tau\\
\\
\hspace{2.8cm}+&\int_O p(t,1,\hat{x}_1;t_0,y)\psi_0(y)dy.
\end{array}
\end{equation}
The function $\phi$ can be represented in terms of the fundamental solution. Define
\begin{equation}
\frac{1}{2}Z^1_F(t,x;\tau ,y)=\frac{\partial}{\partial x_1}p_F(t,x;t_0,y)- p_F(t,x;t_0,y),
\end{equation}
and for $r\geq 1$ define
\begin{equation}
Z^{r+1}_F(t,x,\tau,y)=\int_{t_0}^ t \int_{H_0} Z^ r_F(t,x,\sigma,1,\hat{z}_1)Z^ r_F(\sigma,\hat{z}_1,\tau,y)dH_yd\sigma .
\end{equation}
Then
\begin{equation}
\frac{1}{2}\phi (t,x)=\Gamma(t,x)+\sum_{r=1}^ {\infty} \int_{t_0}^t \int_{H_0} Z^r_F(t,x,\tau,1,\hat{y}_1)dH_y d\tau .
\end{equation}
We summarize
\begin{thm}\label{integrorep} Assume that (GG) and the assumptions of the general framework hold.
Then the free boundary surface function $F$ of the multivariate American basket put option (with weights normalized to $1$ w.l.o.g.) in the transformed coordinates ($t_0\geq 0$)
\begin{equation}
\begin{array}{ll}
\psi: (t_0,T)\times {\mathbb R}^ n_+\rightarrow (t_0,T)\times[1,\infty)\times \left(0,1 \right)^ {n-1}\\
\\
\psi(t,S_1,\cdots ,S_n)= \left( t ,\frac{\sum_{i=1}^n S_i}{F},\frac{S_2}{\sum_{i=1}^n S_i},\cdots ,\frac{S_n}{\sum_{i=1}^n S_i}\right)
\end{array}
\end{equation}
is solution of the integral equation
\begin{equation}\label{freesurf}
\begin{array}{ll}
F(t,x_2,\cdots ,x_n)=K-u(t,1,\hat{x}_1)=\\
\\
\hspace{2cm} K-\int_{t_0}^ t \int_{H_0} p_F(t,1,\hat{x}_1;\tau,y)\phi (\tau,1,\hat{y}_1)dH_{y}d\tau\\
\\
\hspace{2cm}+\int_O p_F(t,1,\hat{x}_1;t_0,y)\psi_0(y)dy,
\end{array}
\end{equation}
where
\begin{equation}
\frac{1}{2}\phi (t,x)=\Gamma(t,x)+\sum_{r=1}^ {\infty} \int_{t_0}^t \int_{H_0} Z^r_F(t,x,\tau,1,\hat{y}_1)dH_y d\tau
\end{equation}
with 
\begin{equation}
\frac{1}{2}Z^1_F(t,x;\tau ,y)=\frac{\partial}{\partial x_1}p_F(t,x;t_0,y)- p_F(t,x;t_0,y),
\end{equation}
\begin{equation}
Z^{r+1}_F(t,x,\tau,y)=\int_0^ t \int_{H_0} Z^ r_F(t,x,\sigma,z)Z^ r_F(\sigma,z,\tau,1,\hat{y}_1)dH_yd\sigma,
\end{equation}
for $r\geq 1$, and $p_F$ is the fundamental solution of
\begin{equation}\label{parabolicF}
u_{t}=\frac{F_t}{F}x_1\frac{\partial u}{\partial x_1}+\frac{1}{2}\sum_{ij}a^F_{ij}\frac{\partial^ 2 u}{\partial x_i\partial x_j}+\sum_j b^F_j\frac{\partial u}{\partial x_j}+rx_1\frac{\partial u}{\partial x_1}.
\end{equation}

\end{thm}

\begin{rem}
The construction of the fundamental solution $p_F$ will be considered in Section 5 and Section 6 below.
Note that the fundamental solution can be approximated accuratively by WKB-expansions. This will be exploited when we derive a scheme from Theorem 4.1. below.
 shall use it for proving existence and uniqueness of the free boundary function below.
\end{rem}

\section{Analysis of the nonlinear equation 3.5}
From the relation 
\begin{equation}\label{freeb}
F(t,\hat{x}_1)=K-u(t,1,\hat{x}_1)
\end{equation}
(recall that $(x_2,\cdots ,x_n)=\hat{x}_1$) it is clear that  uniqueness of $u$ implies uniqueness of $F$ and vice versa. Moreover, the initial value (free)-boundary problem \eqref{frontbd} can be written as a nonlinear initial value boundary problem by substitution of $F$ by $K-u(t,1,\hat{x}_1)$ in the coefficients $a_{ij}^ F$ and $b_i^ F$ and where (BC3) can be dropped. Well, there are powerful techniques in order to show that a weak H\"{o}lder continuous solution exists. Our interest is to prove regularity in a constructive way in order to get a computable scheme. In order to define such a regular scheme for the global solution and the free boundary it is convenient to consider the equation for an equivalent function $v=tu$. From (3.5) we get
 \begin{equation}
\begin{array}{ll}
v_{t}=\frac{F_t}{F}x_1 \frac{\partial v}{\partial x_1}+\frac{1}{2}\sum_{ij}a^F_{ij}\frac{\partial^ 2 v}{\partial x_i\partial x_j}+\sum_j b^F_j\frac{\partial v}{\partial x_j}\\
\\
+r\left( x_1\frac{\partial v}{\partial x_1}-v\right)+u.
\end{array}
\end{equation}
with initial condition $v(0,x)=0$. Note that the original function $u$ appears now as a source term on the right side. The associated boundary conditions (vBC1), (vBC2) and (vBC3) are also obtained by multiplication of (BC1), (BC2), and (BC3) with $t$ followed by substitution of $tu$ by $v$. Especially, (vBC2) becomes 
\begin{equation}
(vBC2)~~v_{x_1}(t,1,x_2,\cdots ,x_n)-v(t,1,x_2,\cdots ,x_n)=-Kt,
\end{equation}
and this matches well with $v(0,.)\equiv 0$ at $t=0$.
We assume that \eqref{frontbd} has been transformed to a finite domain (a suitable additional transformation with respect to the variable $x_1$ s easily constructed) in the sense that $H_{\infty}$ is transformed to $H_d$ with $d>1$. 
Therefore, in the following we shall be able to use some results on bounded domains $D_b\subset {\mathbb R}^{n+1}$ where $D_b$ is the bounded domain obtained after transformation with respect to the first variable. 
We abbreviate $\frac{\partial}{\partial x_i}=:D_{x_i},~\frac{\partial^2}{\partial x_i\partial x_j}=:D_{x_ix_j},~\frac{\partial}{\partial t}=:D_{t}$.
Note that \eqref{frontbd} is equivalent to a problem on a bounded domain with Dirichlet conditions on a transformed hyperplane $H_{\infty}$. Including infinity into the domain we denote the boundary of \eqref{frontbd} by $\partial D=H_1 \cup H_{\infty} \cup I \cup f\cup S_{\mbox{ang}}$
and, for convenience, we define $S:=H_1 \cup H_{\infty}\cup S_{\mbox{ang}}$
with $H_1:=H=\left\lbrace (t,x)\in (0,T)\times \overline{O}|x_1=1\right\rbrace$, \newline $H_{\infty}:=\left\lbrace (t,x)\in (0,T)\times \overline{O}|x_1=\infty\right\rbrace$, $I=\left\lbrace (t,x)\in [0,T]\times \overline{O}|x_1\geq 1, t=0\right\rbrace$,
$f=\left\lbrace (t,x)\in [0,T]\times \overline{O}|x_1\geq 1 , t=T\right\rbrace$, and \newline
$S_{\mbox{ang}}=\left\lbrace (t,x)\in [0,T]\times \overline{O}|x_i\in\left\lbrace 0,1\right\rbrace , i\geq 2\right\rbrace.$
For any points $P\in S$ and $Q\in D_b$ define the distance from $Q$ to $S$ by
\begin{equation}
d_Q:=\inf_{P\in S} d(P,Q),
\end{equation}
where $d(P,Q)$ denotes the Euclidean distance between $P=(t,x)$ and $Q=(s,y)$, i.e. $d(P,Q)=\sqrt{ |x-y|^2+|t-s|^2}.$
Furthermore, for any two points $P,Q\in D_b$, integers $m\geq 1$ and any function $v:D_b\rightarrow {\mathbb R}$ let $d_{PQ}=\min{\left\lbrace d_p,d_Q\right\rbrace  }$ and define 
\begin{equation}
\delta^{\alpha+m}_{\alpha}(d^m v):=\mbox{lub}_{P,Q\in D_b}d_{PQ}^{m+\alpha}\frac{|v(P)-v(Q)|}{d(P,Q)^{\alpha}},
\end{equation}
and
\begin{equation}
|d^mv|^{o}_{\alpha}=\sup_{(t,x)\in D_b}|d^mv(t,x)|+\delta^{m+\alpha}_{\alpha}(d^mv)
\end{equation}
Let
\begin{equation}
|v|^{o}_{2+\alpha}=|v|^{o}_{\alpha}+\sum_{i=1}^n|d D_{x_i}v|^{o}_{\alpha}+\sum_{ij=1}^n|d^ 2D_{x_ix_j}v|^{o}_{\alpha}+
|d^2 D_t v|^{o}_{\alpha},
\end{equation}
and define
\begin{equation}
C^o_{2+\alpha}:=\left\lbrace f:D\rightarrow {\mathbb R}| |f|^{o}_{2+\alpha}<\infty\right\rbrace 
\end{equation}
We also use Banach spaces which are regular on some part of the boundary. Therefore we define
\begin{equation}
|v|_{\alpha}=\sup_{(t,x)\in D_b}|v(t,x)|+\delta^{0}_{\alpha}(v)
\end{equation}
\begin{equation}
|v|_{2+\alpha}=|v|_{\alpha}+\sum_{i=1}^n| D_{x_i}v|_{\alpha}+\sum_{ij=1}^n|D_{x_ix_j}v|_{\alpha}+
|D_t v|_{\alpha},
\end{equation}
and define
\begin{equation}
C_{2+\alpha}:=\left\lbrace f:D\rightarrow {\mathbb R}| |f|_{2+\alpha}<\infty\right\rbrace .
\end{equation}
Note that we keep the domain implicit if this is clear from the context. If we specify some subdomains, e.g.  the suddomain of the angle variables $\left]0,1\right[^{n-1}$, then $C_{2+\alpha}\left(\left]0,1\right[^{n-1} \right) $ refers to the set of functions where spatial derivatives up to second order are H\"{o}lder in the sense of the norm $|w|_{2+\alpha}=|w|_{\alpha}+\sum_{i=1}^n| D_{x_i}w|_{\alpha}+\sum_{ij=1}^n|D_{x_ix_j}w|_{\alpha}$. Analogously, $C^o_{1+\alpha}\left(\left]0,T\right[\right)$ refers to the Banach space with the norm   $|v|^{o}_{1+\alpha}=|v|^{o}_{\alpha}+|d^2 D_t v|^{o}_{\alpha}$ etc.

Recall that
\begin{prop}
If $D_b\subset {\mathbb R}^n$ is bounded, then $C_{2+\alpha}(D_b)$ is a Banach space.
\end{prop}
Next we recall some classical results on a priori estimates by Schauder.
Consider the first initial-boundary value problem on a hypercube $D$: 

\begin{equation}\label{initbound}
\left\lbrace \begin{array}{ll} \frac{\partial w}{\partial  t}-\frac{1}{2}\sum_{ij}a_{ij}(t,x)\frac{\partial^2 w}{\partial x_i\partial x_j}-\sum_i b_i(t,x)\frac{\partial w}{\partial x_i}+c(t,x)w=f(t,x) \\
\hspace{8.5cm}\mbox{ in } D_b,\\
\\
w(t,x)=h(t,x) \mbox{ on } I\cup H\times (0,T)\cup H_d\times (0,T),
\end{array}\right.
\end{equation}
where $H=\left\lbrace x\in \overline{D_b}| x_1=1  \right\rbrace $, $H_d=\left\lbrace x\in \overline{D_b}| x_1=d  \right\rbrace $ (for some $d>1$, and $I$ and $F$ are the boundaries of $D$ at $t=0$ and $t=T$ respectively.
Assume that the following assumptions are satisfied:
\begin{itemize}
\item[(A)] the coefficients $(t,x)\rightarrow a_{ij}(t,x)$ satisfy an ellipticity condition, i.e. there exists a constant $C>0$ such that for any $(t,x)\in D_b$
\begin{equation}
\sum_{ij}a_{ij}(t,x)\xi_i\xi_j\geq c|\xi|^2. 
\end{equation}

\item[(B)] The coefficient functions $(t,x)\rightarrow a_{ij}(t,x)$ and $(t,x)\rightarrow b_i(t,x)$ and $(t,x)\rightarrow c(t,x)$ are locally H\"older continuous (exponent $\alpha$), i.e. there exists a constant $C>0$ such that 
\begin{equation}
|a_{ij}|_{\alpha}\leq C, |db_i|_{\alpha}\leq C, |d^2c|_{\alpha}\leq C
\end{equation}

\item[(C)] the function $f$ is locally H\"older continuous (exponent $\alpha$), i.e. there exists a constant $C>0$ such that 
\begin{equation}
|d^2f|_{\alpha}\leq C
\end{equation}
\end{itemize}

\begin{rem}
We say that (A), and (B) of the present section hold for the domain $D$ for the problem (3.5) if (A) and (B) hold for the equivalent problem defined on a finite domain $D_b$ where $H_{\infty}$ has been transformed to $H_d$.

\end{rem}

Now a classical result on Schauder estimates states that
\begin{thm}
Assume that (A), (B), and (C) hold and assume that $S$ has the outside strong sphere property. Then for any continuous function $h$ on $I\cup H$ there exists a unique solution of the \eqref{initbound}, and $u\in C_{2+\alpha}(D_b)$. 
\end{thm}

\begin{thm}
Assume that $(A), (B), \mbox{ and } (C)$ are satisfied. Then there exists a constant $K$ depending only on the ellipticity constant $c$ and the H\"{o}lder constant $C$ such that any solution of the first equation of (\ref{initbound}) satisfies
\begin{equation}\label{aprioripara}
|w|^o_{1+\delta/2,2+\delta}\leq K_0\left[|w|^o_{0}+|d^2f|^o_{\delta}\right]. 
\end{equation}
\end{thm}
We are going to construct the free boundary $F$ from the fixed point of an iterative cnstruction which leads to the solution of the equation (3.5) in a subspace of $C_{2+\alpha}$. This subspace is
\begin{equation}
B_{2+\alpha}:=C^o_{1+\alpha}\left( \left]0,T\right[ \right)\times C_{\alpha}\left[1,\infty \right]\times   C_{2+\alpha}\left(\left]0,1 \right[^{n-1}  \right)  
\end{equation}
which induces a space on the hyperplane $H$ of form
\begin{equation}
B^H_{2+\alpha}:=C^o_{1+\alpha}\left( \left]0,T\right[ \right)\times   C_{2+\alpha}\left(\left]0,1 \right[^{n-1}  \right)  
\end{equation} 
Now we can establish the main result.
\begin{thm}\label{mainthm}Assume that the assumption of the general framework, the assumptions (A), and (B) of the present section hold for the domain $D$ (cf. remark 5.2) , and the global graph condition (GG) holds.
Then free boundary function $F:]0,T[\times ]0,1[^{n-1}\rightarrow {\mathbb R}$ in \eqref{frontbd} is in $C^o_{2+\alpha} $. More precisely $F$ is in the subspace $B^H_{2+\alpha}\subset C^o_{2+\alpha}$ (here, all Banach spaces are considered with respect to the domain of $F$).
\end{thm}
Proof. Note that only the inner regularity of the free boundary function can be proved. 
Starting with some $u^0$ (solution of (3.5) for $F\equiv1$ for example) we consider for $n\geq 1$ the following iteration for  \ref{frontbd}. Let
\begin{equation}\label{redpara*}
\begin{array}{ll}
v^n_{t}=\frac{F^{n-1}_t}{F^{n-1}}x_1 \frac{\partial v^n}{\partial x_1}+\frac{1}{2}\sum_{ij}a^{F^{n-1}}_{ij}\frac{\partial^ 2 v^n}{\partial x_i\partial x_j}+\sum_j b^{F^{n-1}}_j\frac{\partial v^n}{\partial x_j}\\
\\
+r\left( x_1\frac{\partial v^n}{\partial x_1}-v^n\right)+u^{n-1},
\end{array}
\end{equation}
where $v^n$ has zero initial condition. On the hyperplane $\left\lbrace x_1=1\right\rbrace $ we require
\begin{equation}
v^n_{x_1}(t,1,x_2,\cdots ,x_n)-v^n(t,1,x_2,\cdots ,x_n)=-t K.
\end{equation}
Given $u^{n-1}$ and $F^{n-1}$ for each $n$ we look first at the solution for $v^n$ in the form
(note that $\int_{D_b}$ below contains an integral $\int_0^t$)
\begin{equation}\label{itequation}
\begin{array}{ll}
v^n(t,x)=&\int_0^ t \int_H p_{v^n}(t,1, x;\tau,1, \hat{y}_1)\phi_{v^n} (\tau,1,\hat{y}_1)dH_yd\tau\\
\\
&+\int_{D_b} u^{n-1}(s,y)p_{v^n}(t,x;s,y)dyds,
\end{array}
\end{equation}
where $p_{v^n}$ is the fundamental solution of
\begin{equation}
v^n_{t}=\frac{1}{2}\sum_{ij}a^{F^{n-1}}_{ij}\frac{\partial^ 2 v^n}{\partial x_i\partial x_j}+\sum_j b^{F^{n-1}}_j\frac{\partial v^n}{\partial x_j}+\frac{F^{n-1}_t}{F^{n-1}}x_1 \frac{\partial v^n}{\partial x_1},
\end{equation}
and where $\phi_{v^n}$ solves an integral equation 
\begin{equation}\label{vol2}
\begin{array}{ll}
\frac{1}{2}\phi_{v_n}(t,1,\hat{x}_1)=\Gamma_{v_n}(t,1,\hat{x}_1)+\\
\\ \int_{t_0}^ t \int_{H_0} (\frac{\partial}{\partial x_1}p_{v_n}(t,1,\hat{x}_1,\tau,\hat{y}_1)-p_{v_n}(t,1,\hat{x}_1;\tau,1,\hat{y}_1))\phi_{v_n} (\tau,1,\hat{y}_1)dH_yd\tau,
\end{array}
\end{equation}
(defining $\phi_{v_n}$) with 
\begin{equation}
\begin{array}{rrr}
\Gamma_{v_n} (t,x)=&\int_{O} \left( \frac{\partial}{\partial x_1}p_{v_n}(t,x;t_0,y)- p_{v_n}(t,x;t_0,y)\right) t\psi_0(y)dy+t K.
\end{array}
\end{equation}
\begin{rem}
Classical theory shows that the fundamental solution (or density) $p_F$ exists. This can be also shown by use of WKB-expansions.
\end{rem}
For the next step we get $F^n$ and $u^n$ via $tF_n=K-v^n(1,\hat{x}_1)$ and $tu^n=v^n$.  Next we prove that equation (\ref{itequation}) leads to fixed point equation $G$ for the free boundary.
is used. The next lemma shows that such a fixed point exists and is indeed located in (a subspace of) the function space $C^o_{2+\alpha}$ (with repect to the domain of the free boundary).
\begin{lem}
Equation (\ref{itequation}) induces a map  $G^H: B^H_{2+\alpha}\rightarrow B^H_{2+\alpha}$ which has fixed point defining the free boundary $F$ in $B^H_{2+\alpha}$.
\end{lem}
\begin{proof}
The two terms on the right side of (\ref{itequation}) are equivalent to solutions of Cauchy problems with H\"{o}lder source term and zero initial data. Indeed, the first term on the right side represents the solution of a Cauchy problem on the hyperplane $H=\left\lbrace  x_1=1\right\rbrace $ of the form
\begin{equation}\label{itsecondbd1}
\left\lbrace \begin{array}{ll}
v^n_{t}=\frac{1}{2}\sum_{ij}a^{F^{n-1}}_{ij}\frac{\partial^ 2 v^n}{\partial x_i\partial x_j}+\sum_j b^{F^{n-1}}_j\frac{\partial v^n}{\partial x_j}+\frac{F^{n-1}_t}{F^{n-1}}x_1 \frac{\partial v^n}{\partial x_1}+\phi_{v^n} (\tau,1,\hat{y}_1),~\mbox{on}~H,\\
\\
(IC_n)~~v_n(0,1, \hat{x}_1)=0.
\end{array}\right.
\end{equation}
Note that $\phi_n$ is itself the solution of an integral equation where the representation of the solution involves normal derivatives of $p_{v_n}$, i.e. derivatives with respect to $x_1$. Note that $\phi_{v^n}$ is H\"{o}lder continuous by standard arguments.
\begin{rem}
Note that we have atmost $C^1$ regularity at $x_1=1$ (higher derivatives would involve third derivatives of $p_{v_n}$ in the representation of $\phi_n$ which do not exist. However, in order to prove the regularity of the free boundary we note that we do not need to prove the $\phi_n$ has a first derivative with respect to $x_1$ (corresponding to a second derivative of $p_{v_n}$ with respect to $x_1$ in the representation of the solution for $\phi_n$).   
\end{rem}
Let $v^{n,1}$ be the solution of (\ref{itsecondbd1}).
We get the estimate
\begin{equation}\label{aprioripara1}
|v^{n,1}|_{2+\alpha}\leq K_0|d^2\phi_{v^n}|_{\alpha}. 
\end{equation}
The second term on the right side of  (\ref{itequation}) is equivalent to the solution of the initial value problem
\begin{equation}\label{itsecondbd2}
\left\lbrace \begin{array}{ll}
v^n_{t}=\frac{1}{2}\sum_{ij}a^{F^{n-1}}_{ij}\frac{\partial^ 2 v^n}{\partial x_i\partial x_j}+\sum_j b^{F^{n-1}}_j\frac{\partial v^n}{\partial x_j}+\frac{F^{n-1}_t}{F^{n-1}}x_1 \frac{\partial v^n}{\partial x_1}+u^{n-1},~\mbox{on}~D,\\
\\
(IC_n)~~v_n(0,x)=0.
\end{array}\right.
\end{equation}
Let $v^{n,2}$ be the solution of the problem (\ref{itsecondbd2}).
We get the estimate
\begin{equation}\label{aprioripara2}
|v^{n,2}|_{1+\delta/2,2+\delta}\leq K_0|d^2u^{n-1}|_{\delta}. 
\end{equation}
Hence $u^n$ with $tu^n=v^n$ and $F^n(t,.)=K-u^n(t,1,.)$ are in $C_{2+\alpha}$ and $F^n_t$ is in $C^o_{\alpha}$, where $tF^n_t$ in $C_{\alpha}$. 
\begin{rem}
Classical theory shows that inductively for each $n$ we have $F^n(t,.)\leq \frac{C}{t^{\delta}}$ for some constant $C$ and $\delta\in (1/2,1)$.  
\end{rem}

Hence we get inductively that
$u^n,v^n,F^n$ are in the subspace $B_{2+\alpha}$ for all $n$. A suitable definition of the $G:B_{2+\alpha}\rightarrow B_{2+\alpha}$ is then the map which assigns $v^{n-1}=v^{n-1,1}+v^{n-1,2}$ to $v^{n}=v^{n,1}+v^{n,2}$. The map $G$ induces then the map $G^H$ on $B^H_{2+\alpha}$ (restiction to $\left\lbrace x_1=1\right\rbrace$.
Since the constant $K_0$ depends only on the ellipticity constant $c$ and the H\"{o}lder constant $C$ the existence of a fixed point is straightforward (consider weighted norms and derive exponential convergence). However, from a numerical point of view it makes sense to proceed with a certain time discretization. We may consider a time transform
\begin{equation}\label{transtime}
\begin{array}{ll}
t: [0,\infty)\rightarrow [0,\infty),\\
\\
t(\tau)=\rho \tau .
\end{array}
\end{equation}
Then we get an equation in $\tau$ equivalent to (\ref{redpara*}) where the coefficients of the symbol of the operator become small if $\rho$ is small. We have $\frac{dt}{d\tau}=\rho.$ For the transformed functions $w^{\rho}, u^{\rho}$ with $w^{\rho}(\tau,x)=w(t,x)$ and $u^{\rho}(\tau,x)=u(t,x)$ we get
\begin{equation}
\begin{array}{ll}
v^{\rho}_{\tau}=\frac{1}{2}\rho \frac{F_{\tau}}{F}x_1 \frac{\partial v^{\rho}}{\partial x_1}+\frac{1}{2}\sum_{ij}\rho a^F_{ij}\frac{\partial^ 2 v^{\rho}}{\partial x_i\partial x_j}+\sum_j \rho b^F_j\frac{\partial w^{\rho}}{\partial x_j}\\
\\
+r\rho\left( x_1\frac{\partial v^{\rho}}{\partial x_1}-v^{\rho}\right)+\rho\tau u^{\rho}.
\end{array}
\end{equation}
and apply the a priori estimate for small $\rho$, and then iterate the argument in time. A depper numerical shows that the time step size may be increased as time goes by.
\end{proof}
\begin{rem}
 Note that the estmates for $v^n_t$ imply that $tF^n_t$ is always H\"{o}lder up to the boundary. We could impose this additional restriction on $B^H_{2+\alpha}$ (for $\alpha$ according to remark 5.9 above). 
\end{rem}

\section{Iterative schemes for the free boundary surface and the Greeks}

The contraction map $G$ (or $G^H$) in the regularity proof above leads to various versions of iterative schemes for the free boundary surface $F$, its time derivative and its spatial derivatives up to second order. We describe an algorithmic scheme which has several possible realizations (probabilistic and PDE-schemes) of its subproblems. We shall describe the probabilistic scheme in more detail. PDE-schemes and issues of implementation will be considered elsewhere.

\subsection{Splitting scheme (splitting the mixed boundary problem at each iteration step)}
In the proof of theorem \eqref{mainthm} we observed that the free boundary function $F$ is given by relation (\ref{freeb}) where $u(t,1,\hat{x}_1)$ is computed along with $v(t,1,\hat{x}_1)$. In order to make the computation more stable we may consider an analogous iteration for the function $w=\frac{1}{2}t^2u$ based on the fixed point equation
\begin{equation}\label{basicint}
\begin{array}{ll}
w(t,1,\hat{x}_1)=&\int_0^ t \int_H p_w(t,1,\hat{x}_1;\tau,1,\hat{y}_1)\phi_w (\tau,1,\hat{y}_1)dH_yd\tau\\
\\
&+\int_D t u(s,y)p_w(t,1,\hat{x}_1;s,y)dyds,
\end{array}
\end{equation}  
and along with a linear integral equation for $\phi_w$. Now given a function $F\in C^H_{2+\alpha}$ (more precisely, in $B^H_{2+\alpha}$) we may first solve for the linear integral equation for $\phi_w$. This is a Volterra integral equation which is well-studied in the literature numerically (cf. \cite{Sa} for a probabilistic treatment).   
Next we analyze the terms on the right side of \eqref{basicint}. First, the term
\begin{equation}\label{part1}
\int_D \frac{1}{2}su(s,y)p_w(t,1,\hat{x}_1;s,y)dyds
\end{equation}
can be interpreted to be the solution of
\begin{equation}\label{itsecondbd}
\left\lbrace \begin{array}{ll}
w_{t}=\frac{1}{2}\sum_{ij}a^F_{ij}\frac{\partial^ 2 w}{\partial x_i\partial x_j}+\sum_j b^F_j\frac{\partial w}{\partial x_j}+rx_1\frac{\partial w}{\partial x_1}-rw+tu(t,y),\\
\\
(IC)~~w(0,x)=0
\end{array}\right.
\end{equation}
on $D$ with natural boundary conditions on $H$. (Note that a transformation of spatial variables of the form $x\rightarrow z$ with $z_1=\ln \ln (x_1)$, and similar transformations with respect to the other coordinates $\hat{x}_1\rightarrow \hat{z}_1$ shows that the problem is equivalent to a Cauchy problem on the whole space.) 
The interpretation of
\begin{equation}\label{itfirst}
\int_0^ t \int_H p_w(t,1,\hat{x}_1;\tau,1,\hat{y}_1)\phi_w (\tau,1,\hat{y}_1)dH_yd\tau
\end{equation}
as a solution of a Cauchy problem with initial condition equal to zero and a source related to $\phi_w$ of reduced dimension $n-1$
is more theoretical in the sense that involves the solution of an integral equation $\phi_w$ which requires the knowledge of the density (and its normal derivative, i.e. the normal derivative with respect to the hypersurface $H$). This is a first reason why an efficient computation of the transition density (fundamental solution) is desirable, and it also indicates that our algorithm is intrinsically probabilistic to some extent in the sense that it cannot be completely realized without referring to the transition density. A second reason is the following.
If we know the fundamental solution $p_w$ in \eqref{part1} and in \eqref{itsecondbd} (for fixed $F$), then we can obtain the expressions
$p_{w,t}$, $p_{w,x_i}, i\in\{2,\cdots ,n\}$, and $p_{w,x_ix_j}, i,j\in\{2,\cdots ,n\}$ for the next iteration step $n+1$ 
by explicit differentiation, i.e  given the approximation of the free boundary surface of the $n$-th iteration step $F^n$ and its derivatives occurring in $p_{w_n}$ we get the approximation of the free boundary surface of the $n+1$-th iteration step $F^{n+1}$ and its derivatives $F^{n+1}_t$, $F^{n+1}_{x_i}, i\in\{2,\cdots ,n\}$, and $F^{n+1}_{x_ix_j}, i,j\in\{2,\cdots ,n\}$
via derivatives of $w_{n+1}(t,1,\hat{x}_1)$ which can be computed via derivatives of $p_{w_n}$.

As we observed in \cite{KKS}, we may approximate derivatives of value functions by derivatives of analytical approximations of transition densities in convolution representations. These approximations of value functions are numerically efficient and allow for error estimates in strong norms by using a priori estimates of Safonov type (cf. \cite{Kryl}, \cite{KKS}). The WKB-expansion of the transition density is such an analytical representation.

\subsection{Use of WKB-expansions of the fundamental solution at each iteration step}

We review some recent research on the fundamental solution (transition density), i.e.
some results concerning WKB-expansions of parabolic equations (cf. \cite{Ka2, KKS}, for more  details). 
Consider the solution $(t,x,s,y)\rightarrow p(t,x,s,y)$ of the family of parabolic equations
\begin{equation}\label{parabeq}
\begin{array}{l}
\frac{\partial u}{\partial t}-\frac{1}{2}\sum_{i,j}a_{ij}\frac{\partial^2 u}{\partial x_i\partial x_j}-
\sum_i b_i\frac{\partial u}{\partial x_i}=0,\\
\\
p(0,x,0,y)=\delta(x-y) ,
\end{array}
\end{equation}
on a domain $D$, parameterized by $y\in {\mathbb R}^n$, and where $\delta$ denotes the Dirac delta distribution 
 and the diffusion coefficients $a_{ij}$ and the first order coefficients $b_i$ in (\ref{parabeq})
may depend on time $t$ and the spatial variable $x$. Without loss of generality and for simplicity of notation we consider the case where the coefficients depend on the spatial coordinates. Note that the time coordinate in the coefficients can be treated as an extra spatial coordinate and the resulting degenerate parabolic equation belongs to a class of so-called projective parabolic equations which are subject to all the following results (cf. \cite{Ka3}), In the following let $\delta t:=t-s$, and let the functions
$$
(x,y)\rightarrow d(x,y)\ge0,~~(x,y)\rightarrow c_k(x,y),~k\geq 0,
$$
be defined on  $[0,T]\times{\mathbb R}^n\times [0,T]\times{\mathbb R}^n$. 
Then a set of (simplified) conditions sufficient for pointwise valid local WKB-representations 
of the form
\begin{equation}\label{WKBrep}
p(t,x,s,y)=\frac{1}{\sqrt{2\pi \delta t}^n}\exp\left(-\frac{d^2(x,y)}{2\delta t}+\sum_{k= 0}^{\infty}c_k(x,y)\delta t^k\right), 
\end{equation}
is given by
\begin{itemize}
\item[(WKB1)] The operator $L$ is uniformly elliptic in ${\mathbb R}^n$, i.e. the matrix norm of $(a_{ij}(t,x))$ is bounded from below 
by $\lambda>0$ and from above by $\Lambda>\lambda,$ uniformly in~$x$,

\item[(WKB2)] the smooth functions $(t,x)\rightarrow a_{ij}(t,x)$ and $(t,x)\rightarrow b_i(t,x)$ and all their derivatives are bounded.

\end{itemize}

Summing up we have the following theorem:
\begin{thm}
If the hypotheses (A),(B) are satisfied, then the fundamental solution $p$ has the local representation
\begin{equation}\label{repwkb}
p(\delta t,x,y)=\frac{1}{\sqrt{2\pi \delta t}^n}\exp\left( -\frac{d^2(x,y)}{2\delta t}+\sum_{k\geq 0}c_k(x,y)\delta t^k\right),
\end{equation}
where $d$ and $c_k$ are smooth functions, which are unique global solutions of the first order differential equations (\ref{ed}),(\ref{c01e}), and (\ref{1gaa}) below. Especially,
$$(\delta t,x,y)\rightarrow \delta t \ln p(\delta t,x,y) =-\frac{n}{2} \delta t\ln (2\pi \delta t) -\frac{d^2}{2} +\sum_{k\geq 0}c_k(x,y)\delta t^{k+1}$$
is a smooth function which converges to $-\frac{d^2}{2}$ as $\delta t\searrow 0$, where $d$ is the Riemannian distance induced by the line element $ds^2=\sum_{ij}a^{-1}_{ij}dx_idx_j$, where with a slight abuse of notation $(a^{-1}_{ij})$ denotes the matrix inverse of $(a_{ij})$.
\end{thm}
The recursion formulas for $d$ and $c_k,~k\geq 0$ are obtained by plugging  the ansatz (\ref{WKBrep}) into the parabolic equation (\ref{parabeq}), 
and ordering terms with respect to the monoms $\delta t^i=(T-t)^i$ for $i\geq -2$.
By collecting terms of order $\delta t^{-2}$ we obtain
\begin{equation}\label{ed}
d^2=\frac{1}{4}\sum_{ij}d^2_{x_i}a_{ij}d^2_{x_j},
\end{equation}
where $d^2_{x_k}$ denotes the derivative of the function $d^2$ with respect to the variable $x_k$,  with the boundary condition  $d(x,y)=0$ for $x=y.$ 
Collecting terms of order $\delta t^{-1}$ yields
\begin{equation}\label{c01e}
-\frac{n}{2}+\frac{1}{2}Ld^2+\frac{1}{2}\sum_{i} \left( \sum_j\left( a_{ij}(x)+a_{ji}(x)\right) \frac{d^2_{x_j}}{2}\right) \frac{\partial c_{0}}{\partial x_i}(x,y)=0,
\end{equation}
where the boundary condition 
\begin{equation}\label{c01b}
c_0(y,y)=-\frac{1}{2}\ln \sqrt{\mbox{det}\left(a_{ij}(y) \right) }
\end{equation}
determines $c_0$ uniquely for each $y\in {\mathbb R}^n$. Finally,
for $k+1\geq 1$ we obtain
\begin{equation}\label{1gaa}
\begin{array}{ll}
(k+1)c_{k+1}(x,y)+\frac{1}{2}\sum_{ij} a_{ij}(x)\Big(
\frac{d^2_{x_i}}{2}\frac{\partial c_{k+1}}{\partial x_j}
+\frac{d^2_{x_j}}{2} \frac{\partial c_{k+1}}{\partial x_i}\Big)\\
\\
=\frac{1}{2}\sum_{ij}a_{ij}(x)\sum_{l=0}^{k}\frac{\partial c_l}{\partial x_i} \frac{\partial c_{k-l}}{\partial x_j}
+\frac{1}{2}\sum_{ij}a_{ij}(x)\frac{\partial^2 c_k}{\partial x_i\partial x_j}  
+\sum_i b_i(x)\frac{\partial c_{k}}{\partial x_i},
\end{array}
\end{equation}
with boundary conditions
\begin{equation}\label{Rk}
c_{k+1}(x,y)=R_k(y,y) \mbox{ if }~~x=y,
\end{equation}
$R_k$ being the right side of (\ref{1gaa}). In \cite{Ka3} it is shown how the function $d^2$ can be approximated in regular norms if only (WKB1) and (WKB2) are satisfied. The technique can be combined with the regular polynomial interpolations developed in \cite{Ka1}.

\subsection{Description of algorithm (including to sparse grids and weighted Monte-Carlo versions)}

We can solve \eqref{frontbd} by an iterative numerical procedure.

We assume that time steps $t_i$ are small enough such that the map $G^H$ is a contraction. We only consider one time step (simply denoted by $t$), since the procedure is the same for all time steps.

\begin{itemize}
\item[(Step1)] Solve
\begin{equation}\label{step1}
\left\lbrace \begin{array}{ll}
u^1_{1,t}=\frac{1}{2}\sum_{ij}a^1_{ij}\frac{\partial^ 2 u^1_1}{\partial x_i\partial x_j}+\sum_j b^1_j\frac{\partial u^1_1}{\partial x_j}+rx_1\frac{\partial u^1_1}{\partial x_1},\\
\\
(IC)~~u^1_1(0,x)=\max \{K-x_1, 0\},
\end{array}\right.
\end{equation}
where $a^1_{ij}:=a^F_{ij}$ $b^1_j:=b^F_j$ for $F\equiv 1$, i.e. compute
\begin{equation}\label{u11}
u^1_1(t,x):=\int_O p_1(t,x;0,y)\max\left\lbrace K-y_1,0\right\rbrace dy,
\end{equation}
where $p_1(t,x;0,y):=p_F(t,x;0,y)$ with $F\equiv 1$.

\vspace{0.5cm}

Next solve for $\phi_1$ in

\begin{equation}\label{vterra0}
\begin{array}{ll}
\frac{1}{2}\phi^1(t,x)=\Gamma(t,x)+\\
\\ \int_0^ t \int_{H_0} (\frac{\partial}{\partial x_1}p_1(t,1,\hat{x}_1,\tau,\hat{y}_1)+p_1(t,1,\hat{x}_1;\tau,1,\hat{y}_1))\times\\
\\
\phi^1 (\tau,1,\hat{y}_1)dH_yd\tau,
\end{array}
\end{equation}
where $\phi^1=\phi^F$ for $F=1$, and $\Gamma$ analogous as in the previous section.
\vspace{0.5cm}

Next compute
\begin{equation}\label{u12}
u^1_2(t,\hat{x}_1):=\int_0^ t \int_{H_0} p_1(t,1,\hat{x}_1;\tau,1,\hat{y}_1)\phi_F (\tau,1,\hat{y}_1)dH_yd\tau
\end{equation}
and
\begin{equation}
\begin{array}{lll}
F^1(t,\hat{x}_1)=K-u^1_2(t,\hat{x}_1)+u^1_1(t,x).
\end{array}
\end{equation}
Next compute the time (first order) and spatial derivatives (up to second order) of $F_1$ by differentiation of the density $p_1$.

\begin{rem}
For higher dimension models the integrals for $u^1_1$ in (\ref{u11}) and $u^1_2$ in (\ref{u12}) have to be computed via MC methods (cf. (\cite{KKS, FrKa, FK, S})  
\end{rem}

\begin{rem}
The use of WKB-expansions of $p_1$ seems the most efficient method. In \cite{KKS} it was observed that a WKB approximation of the transition density including the term $c_1$ allows pricing of options and its derivatives with respect to the underlyings with maturity of $10$ years in LIBOR models of dimension $n=20$ with one time step (beating all concurrent methods). Moreover, efficiency and accuracy is kept for the derivatives because of explicit WKB-approximations.
\end{rem}

\item[(Step2)] Having computed $F^n$, $w_n$ and $u_n$ for $n\geq 1$ compute $w_{n+1}$, $u_{n+1}$, and $F^{n+1}$ as follows. Solve for
\begin{equation}\label{basicintwn1}
\begin{array}{ll}
&w_{n+1}(t,x)=\\
\\
&\int_0^ t \int_H p_{w_n}(t,1,\hat{x}_1;\tau,1,\hat{y}_1)\phi_{w_n} (\tau,1,\hat{y}_1)dH_yd\tau\\
\\
&+\int_D su_n(s,y)p_{w_n}(t,x;s,y)dyds,
\end{array}
\end{equation}
where $\phi_{w_n}$ is the solution of 
\begin{equation}\label{vterrait}
\begin{array}{ll}
\frac{1}{2}\phi_{w_n}(t,x)=\Gamma_n(t,x)+\\
\\ \int_0^ t \int_{H_0} (\frac{\partial}{\partial x_1}p_{w_n}(t,1,\hat{x}_1,\tau,\hat{y}_1)+p_{w_n}(t,1,\hat{x}_1;\tau,1,\hat{y}_1))\times\\
\\
\phi_{w_n} (\tau,1,\hat{y}_1)dH_yd\tau,
\end{array}
\end{equation}
and where $\Gamma_n$ is defined analogously.
Next we have $t^2u_{n+1}(t,x)=w_{n+1}(t,x)$ (note that a priori estimates imply inductively that $u^{n+1}\in C_{\alpha}\cap C^o_{2+\alpha}$).
Then we have
\begin{equation}
F^{n+1}(t,\hat{x}_1)=K-u_{n+1}(t,1,\hat{x}_1)
\end{equation}
Next compute the time (first order) and spatial derivatives (up to second order) of $F^{n+1}$, and proceed to the next step $n+1$.

\begin{rem}
For higher dimension models the integrals for $w^n_1$ in (\ref{basicintwn1})  have to be computed via MC methods (cf. (\cite{KKS, FrKa, FK, S})  
\end{rem}

\item[(Step3)] The functions $F^n$ form a Cauchy sequence in the Banach space $B^H_{2+\alpha}$. Hence we iterate the step 2 until

\begin{equation}
|F^{n+m}-F^n|_{2+\alpha}\leq \frac{k^n}{1-k}|F_{1}-1|_{2+\alpha}\leq \epsilon
\end{equation}
for a prescribed $\epsilon >0$.

\end{itemize}

\begin{rem}
The initial step 1 and the iterated step 2 performed with WKB approximations may be done with different methods based on the dimension of the problem:
\begin{itemize}
\item[M1] For problems up to dimension $3$ computation grids such as $UG$ may be used. 

\item[M2] For higher dimensional models up to dimension $n=5,6$ sparse grid techniques may be used (cf. \cite{PeSch}. At present -at least to my knowledge- nobody came up with stable numerical sparse grids solutions of higher dimension for linear parabolic problems with a comparable complexity.

\item[M3] In any case Monte-Carlo realizations of the algorithmic scheme presented above are possible. Recently weighted Monte-Carlo schemes have been developed in \cite{FrKa} and \cite{KKS},
where the improved estimators established in the latter article allow also to deal with highly peaked densities which occur especially for small time or small volatility.

\end{itemize}

\end{rem}

\begin{rem}
For the implementation of the algorithm some recent results on regular polynomial implementation (cf. \cite{Ka2}) and the computation of the Riemannian metric (cf. \cite{Ka3}) in regular norms are needed. 
\end{rem}

\begin{rem}
Let us look at the case of higher dimension and probabilistic realizations of the algorithmic scheme presented above.
Given a free boundary approximation function $F^n$ the Monte-Carlo estimators used for the computations of higher dimensional integrals related to linear parabolic Cauchy problems are based on the formula 

\begin{equation}
\frac{\partial^{\alpha} I}{\partial x^{\alpha}}(x)=E\,\frac{\partial^{\alpha}}{\partial
x^{\alpha}}\frac{p_{F_n}(t,x,g(x,\xi))u(g(x,\xi))}{\phi(t,x,g(x,\xi))},\label{Istr},%
\end{equation}
and $\frac{\partial^{\alpha}}{\partial x^{\alpha}}$ are spatial derivatives up to secons order or a similar formula for the time derivative. Here 
$\xi$ be an $\mathbb{R}^{n}$-valued random variable on some probability space with a density $\lambda(z)\neq 0$ for all $z$, 
and the regular (at least twice continuously differentiable) map $\zeta^{x}:=g(x,\xi)$ satisfies $\left\vert
\partial g(x,z)/\partial z\right\vert \neq0,$ and has
density $\phi(x,\cdot)$ on $\mathbb{R}_{+}^{n}$.
The corresponding Monte Carlo estimator is
\begin{equation}\label{est}
\widehat{\frac{\partial^{\alpha} I}{\partial x^{\alpha}}}(x)=\frac{1}{M}\sum_{m=1}^{M}%
\frac{\partial}{\partial x}\frac{p(x,g(x,_{m}\xi))u(g(x,_{m}\xi))}%
{\phi(x,g(x,_{m}\xi))}.%
\end{equation}
In \cite{KKS} it is shown that the latter estimator is of bounded variance even for small time or volatility. Moreover, Monte Carlo methods for computing \eqref{vterrait}  can be found in \cite{Sa}. 
\end{rem}

\section{Epilog}

We have proved regularity of the free boundary surface for a considerable class of relevant market models, and we have set up in detail a scheme for the important and difficult problem of computing the Greeks for American type options which works in the context of higher dimension. On the way we have constructed a nonlinear integral equation which characterizes the free boundary in the multivariate case. For the implementation of the scheme in a general situation an implementation of the computation of the Varadhan metric in regular norms and of the drift functions is needed. This is possible using regular polynomial interpolation (cf. \cite{Ka1} and the analysis of the eikonal equation characterizing the Varadhan metric global approximation of its solution (cf. \cite{Ka3}). The scheme established is very flexible. Note that PDE-scheme realizations are possible beside Monte-carlo realizations especially for lower-dimensional models. Detailed error analysis for different realizations of the scheme and their implementations is certainly of interest. Furthermore, extension to some class of jump-diffusion models is possible and will be considered elsewhere.

\end{document}